\documentclass[a4paper,11pt,reqno]{amsart}

\usepackage{lscape}
\usepackage{multirow}
\usepackage[all]{xy}
\usepackage{amsmath,amssymb,graphicx,enumerate,amsthm}

\def\pr{\mathrm{pr}}
\def\tr{\mathrm{tr}}
\def\<#1,#2>{\langle\,#1,\,#2\,\rangle}

\def\dim{\mathrm{dim}}

\def\s#1{\mathrm{\Sigma_{#1} M}}
\def\sigm{{\mathrm{\Sigma M}}}

\def\qed{\ensuremath{\hfill\Box}}

\def\Aut{\mathrm{Aut}}
\def\so{\mathfrak{so}}
\def\SO{\mathrm{SO}}
\def\CO{\mathrm{CO}}

\def\Sp{\mathrm{Sp}}

\def\SU{\mathrm{SU}}

\def\Spin{\mathrm{Spin}}

\def\GL{\mathrm{GL}}

\def\1{\mathbf{1}}
\def\#{\sharp}

\def\Aut{\mathrm{Aut}}

\def\pr{\mathrm{pr}}

\def\tr{\mathrm{tr}}

\def\wb{Weitzenb\"ock }

\def\<#1,#2>{\langle\,#1,\,#2\,\rangle}

\def\rectangle(#1,#2)[#3,#4]#5{
 \multiput(#1,#2)(#3,0)2{\line(0,1){#4}}\multiput(#1,#2)(0,#4)2{\line(1,0){#3}}
 \put(#1,#2){\vbox to #4pt{\hbox to #3pt{\hfill}\vfill}}}
\def\recttext(#1,#2)[#3,#4]#5{\put(#1,#2)
 {\vbox to #4pt{\vfill\hbox to #3pt{\hss#5\hss}\vfill}}}

\def\qed{\ensuremath{\hfill\Box}}
\newtheorem{Lemma}{Lemma}[section]
\newtheorem{Proposition}[Lemma]{Proposition}
\newtheorem{Theorem}[Lemma]{Theorem}
\newtheorem{Corollary}[Lemma]{Corollary}

\theoremstyle{definition}
\newtheorem{Definition}[Lemma]{Definition}
\newtheorem{Example}[Lemma]{Example}

\newtheorem{Remark}[Lemma]{Remark}

\def\vs{\vskip .4cm}

\numberwithin{equation}{section}

\title[A Note on the Conformal Invariance of $G$-Generalized Gradients]{A Note on the Conformal Invariance of\\
 $G$-Generalized Gradients}
\author{Mihaela Pilca}
\thanks{The financial support of the Graduate School 1269 ``Global Structures in Geometry and Analysis" supported by DFG and the Mathematical Institute of University of Cologne is gratefully acknowledged.}
\address{Mihaela Pilca \\ Mathematisches Institut\\ Universit\"at zu K\"oln\\ Weyertal 86-90 D-50931 K\"oln\\ Germany.}
\email{mpilca@mi.uni-koeln.de}

\begin{document}

\begin{abstract}
We consider generalized gradients in the general context of $G$-structures. They are natural first order differential operators acting on sections of vector bundles associated to irreducible $G$-representations. We study their geometric properties and show in particular their conformal invariance.

\vs

\end{abstract}

\maketitle

\section{Introduction}

The purpose of this note is to introduce generalized gradients associated to $G$-structures and to study their geometric properties. Our main result is their conformal invariance (see Propositions~\ref{confweightg} and \ref{equivg}). 

The classical notion of generalized gradients, also called Stein-Weiss operators, was first introduced by Stein and Weiss, \cite{sw}, on an oriented Riemannian manifold, as a generalization of the Cauchy-Riemann equations. They are first order differential operators acting on sections of vector bundles associated to irreducible representations of the special orthogonal group (or of the spin group, if the manifold is spin), which are given by the following universal construction: one projects onto an irreducible subbundle the covariant derivative induced on the associated vector bundle by the Levi-Civita connection (or, more generally, by any metric connection).

Some of the most important first order differential operators which naturally appear in geometry are, up to normalization, particular cases of generalized gradients. For example, on a Riemannian manifold, the exterior differential acting on differential forms, its formal adjoint, the codifferential, and the conformal Killing operator on $1$-forms are generalized gradients. On a spin manifold, classical examples of generalized gradients are the Dirac operator, the twistor (or Penrose) operator and the Rarita-Schwinger operator. 

An essential property of generalized gradients is their invariance at conformal changes of the metric. This property was noticed for the first time in 1974 by Hitchin, \cite{hit74}, in the case of the Dirac and the twistor operator in spin geometry and it turned out to have important consequences in physics. Two years later, Fegan, \cite{feg}, showed that, up to the composition with a bundle map, the only conformally invariant first order differential operators between vector bundles associated to the bundle of oriented frames are the generalized gradients. Further results in this direction were obtained by Homma, \cite{h1}, \cite{h2}, \cite{h4}, for the conformal invariance of generalized gradients associated to $\mathrm{U}(n)$, $\Sp(n)$ and $\Sp(1)\!\cdot\!\Sp(n)$-structures. For these subgroups, Homma's proof of conformal invariance is given by explicit computations based on the relationship between the enveloping algebra of the Lie algebra of $G$ and the algebraic structure of the principal symbols of generalized gradients.

We propose in this note a unitary proof of the conformal invariance for the known cases and, more generally, for all $G$-generalized gradients, where the $G$-structure is any of the interesting geometric structures that are mostly encountered in literature, \emph{i.e.} $G$ is one of the groups: $\SO(n)$, $\mathrm{U}(\frac{n}{2})$, $\SU(\frac{n}{2})$, $\Sp(\frac{n}{4})$, $\Sp(1)\!\cdot\!\Sp(\frac{n}{4})$, $G_2$ or $\Spin(7)$. The main tools of our proof are the Weyl structures and the so-called conformal weight operator. We mention that such an approach was suggested by Gauduchon, \cite{pg1}, for the case of $\SO(n)$-generalized gradients. We also give some consequences and applications of the conformal invariance of $G$-generalized gradients.

\noindent \textsc{acknowledgment.} This note is a part of my Ph.D. thesis. I thank very much my supervisor Uwe Semmelmann for suggesting me this topic.

\section{Generalized Gradients of $G$-Structures}\label{sectgengrad}

We first briefly recall the general construction of generalized gradients given by Stein and Weiss, \cite{sw}, on an oriented Riemannian (spin) manifold, \emph{i.e.} for the structure group $\SO(n)$ (or $\Spin(n)$), and then we show how this construction can be carried over to $G$-structures.

\subsection{$\SO(n)$ and $\Spin(n)$-Generalized Gradients} 
Let us first state the general context, fix the notations and briefly recall the representation theoretical background needed to define the generalized gradients. The description of the representations of $\so(n)$, the Lie algebra of $\SO(n)$, differs slightly according to the parity of $n$. We write $n=2m$ if $n$ is even and $n=2m+1$ if $n$ is odd, where $m$ is the rank of $\so(n)$. Let $\{e_1, \dots,e_n\}$ be a fixed oriented orthonormal basis  of $\mathbb{R}^n$, so that $\{e_i\wedge e_j\}_{i<j}$ is a basis of the Lie algebra $\so(n)\cong\Lambda^2\mathbb{R}^n$. We also fix a Cartan subalgebra $\mathfrak{h}$ of $\so(n)$ by the basis $\{e_1\wedge e_2, \dots, e_{2m-1}\wedge e_{2m}\}$ and denote the dual basis of $\mathfrak{h}^*$ by $\{\varepsilon_1,\dots,\varepsilon_m\}$. The Killing form is normalized such that this basis is orthonormal. Roots and weights are given by their coordinates with respect to the orthonormal basis $\{\varepsilon_i\}_{i=\overline{1,m}}$. Finite-dimensional complex irreducible $\so(n)$-representations are parametrized by \emph{dominant weights}, \emph{i.e.} those weights whose coordinates are either all integers or all half-integers,\linebreak $\lambda=(\lambda_1,\dots,\lambda_m)\in\mathbb{Z}^m\cup (\frac{1}{2}+\mathbb{Z})^m$ and which satisfy the inequality:
\begin{equation}\label{domw}
\begin{split}
\lambda_1\geq \lambda_2 \geq \cdots \lambda_{m-1}\geq |\lambda_m|, \quad &\text{if }n=2m, \text{ or }\\
\lambda_1\geq \lambda_2 \geq \cdots \lambda_{m-1}\geq \lambda_m \geq 0, \quad &\text{if }n=2m+1.
\end{split}
\end{equation}
Through this parametrization a dominant weight $\lambda$ is the highest weight of the corresponding representation. With a slight abuse of notation, we use the same symbol for an irreducible representation and its highest weight. The representations of $\so(n)$ are in one-to-one correspondence with the representations of the corresponding simply-connected Lie group, \emph{i.e.} $\Spin(n)$, the universal covering of $\SO(n)$. The representations which factor through $\SO(n)$ are exactly those with $\lambda\in\mathbb{Z}^m$. For example, the (complex) standard representation, denoted by $\tau$, is given by the weight $(1,0,\dots,0)$; the weight $(1,\dots,1,0,\dots,0)$ (with $p$ ones) corresponds to the $p$-form representation $\Lambda^p\mathbb{R}^n$, whereas the dominant weights $\lambda=(1,\dots,1,\pm1)$, for $n=2m$, correspond to selfdual, respectively antiselfdual $m$-forms; the representation of totally symmetric traceless tensors $S^p_0\mathbb{R}^n$ has highest weight $(p,0,\dots,0)$. 

The following so-called \emph{classical selection rule}\index{classical selection rule} (see \cite{feg}) describes the decomposition of the tensor product $\tau\otimes\lambda$ into irreducible $\so(n)$-representations, where $\tau$ is the standard representation and $\lambda$ is any irreducible representation. 
\begin{Lemma}\label{selectrule}
An irreducible representation of highest weight $\mu$ occurs in the decomposition of $\tau\otimes\lambda$ if and only if the following two conditions are fulfilled:
\begin{itemize}
	\item[(i)] $\mu=\lambda\pm\varepsilon_j$, for some $j=1,\dots, m$, or $n=2m+1$, $\lambda_m>0$ and $\mu=\lambda$,	
  \item[(ii)] $\mu$ is a dominant weight, \emph{i.e.} satisfies the inequality \eqref{domw}.
\end{itemize}
\end{Lemma}
We adopt the same terminology as in \cite{usgw2} and call \emph{relevant weights of $\lambda$} (and write $\varepsilon\subset\lambda$) those weights $\varepsilon$ of $\tau$, $\varepsilon\in\{0,\pm\varepsilon_1,\dots,\pm\varepsilon_m\}$, with the property that $\lambda+\varepsilon$ occurs in the decomposition of $\tau\otimes\lambda$. The decomposition of the tensor product is then expressed as follows:
\begin{equation}\label{decompirred}
\tau\otimes\lambda=\underset{\varepsilon\subset\lambda}{\oplus}(\lambda+\varepsilon).
\end{equation}
The essential property of the decomposition \eqref{decompirred} is that it is multiplicity-free, \emph{i.e.} the isotypical components are actually irreducible. It thus follows that the projections onto each irreducible summand $\lambda+\varepsilon$ in the splitting are well-defined; we denote them by $\Pi_\varepsilon$.

Let $(M,g)$ be a Riemannian manifold, $\SO_g M$ denotes the principal $\SO(n)$-bundle of oriented orthonormal frames and $\nabla$ any metric connection, considered either as a connection $1$-form on $\SO_g M$ or as a covariant derivative on the tangent bundle $\mathrm{TM}$. If $M$ has, in addition, a spin structure, then we consider the corresponding principal $\Spin(n)$-bundle, denoted by $\Spin_g M$, and the induced metric connection. We consider vector bundles associated to $\SO_g M$ (or $\Spin_g M$) and irreducible $\SO(n)$ ($\Spin(n)$)-representations of highest weight $\lambda$ and denote them by $V_\lambda M$. For instance, the tangent bundle is associated to the standard representation $\tau:\SO(n)\hookrightarrow \GL(\mathbb{R}^n)$ and the bundle of $p$-forms is associated to the irreducible representation of highest weight $\lambda=(1,\dots,1,0,\dots,0)$ (with $p$ ones). We identify the cotangent bundle $\mathrm{T}^*\mathrm{M}$ and the tangent bundle $\mathrm{TM}$ using the metric $g$, since they are associated to equivalent $\SO(n)$-representations. The decomposition \eqref{decompirred} carries over to the associated vector bundles:
\begin{equation}\label{decirrmfd}
\mathrm{T}^*\mathrm{M}\otimes V_\lambda M\cong \mathrm{TM}\otimes V_\lambda M\cong \underset{\varepsilon\subset\lambda}{\oplus}V_{\lambda+\varepsilon} M
\end{equation}
and the corresponding projections are also denoted by $\Pi_\varepsilon$.

A metric connection $\nabla$ on $\SO_g M$ (or $\Spin_g M$) induces on any associated vector bundle $V_\lambda M$ a covariant derivative, denoted by $\nabla^\lambda:\Gamma(V_\lambda M)\to\Gamma(\mathrm{T}^*\mathrm{M}\otimes V_\lambda M)$. The generalized gradients are then built-up by projecting $\nabla^\lambda$  onto the irreducible subbundles $V_{\lambda+\varepsilon} M$ given by the splitting \eqref{decirrmfd}. 
\begin{Definition}\label{defingg}
Let $(M,g)$ be a Riemannian manifold, $\nabla$ a metric connection and $V_\lambda M$ the vector bundle associated to the irreducible $\SO(n)$ (or $\Spin(n)$)-representation of highest \mbox{weight $\lambda$}. For each relevant weight $\varepsilon$ of $\lambda$, \emph{i.e.} for each irreducible component in the decomposition of $\mathrm{T}^*\mathrm{M}\otimes V_\lambda M$, there is a \emph{generalized gradient} $P_\varepsilon^{\nabla^\lambda}$ defined by the composition:
\begin{equation}\label{defgg}
\Gamma(V_\lambda M) \xrightarrow{\nabla^\lambda} \Gamma(\mathrm{TM}\otimes V_\lambda M) \xrightarrow{\Pi_\varepsilon} \Gamma(V_{\lambda+\varepsilon} M). 
\end{equation}
\end{Definition}

The classical case is when $\nabla$ is the Levi-Civita connection. The examples given in the sequel are of this type. However, generalized gradients may be defined by any metric connection. Those defined by the Levi-Civita connection play an important role since they are conformal invariant (see \S ~\ref{confinv}).

\begin{Example}[Generalized Gradients on Differential Forms]\label{pforms}
We consider the bundle of $p$-forms, $\Lambda^p M$, on a Riemannian manifold $(M^n,g)$ and assume for simplicity that $n=2m+1$ and $p\leq m-1$. The highest weight of the representation is $\lambda_p=(1,\dots,1,0,\dots,0)$ and, by the selection rule in Lemma~\ref{selectrule}, it follows that there are three relevant weights for $\lambda_p$, namely $-\varepsilon_p$, $\varepsilon_{p+1}$ and $\varepsilon_1$. The tensor product then decomposes as follows:
\[\mathrm{TM} \otimes\Lambda^p M\cong\Lambda^{p-1}M\oplus \Lambda^{p+1} M\oplus \Lambda^{p,1} M,\]
where the last irreducible component is the Cartan summand (whose highest weight is equal to the sum of the highest weights of the factors of the tensor product). The generalized gradients in this case are, up to a constant factor, the following:  the codifferential, $\delta$, the exterior derivative, $d$, and respectively the so-called \emph{twistor operator}, $T$.
\end{Example}

\begin{Example}[Dirac and Twistor Operator]\label{exsp}
The spinor representation \mbox{$\rho_n\!:\! \Spin(n) \to\! \Aut(\Sigma_n)$}, with $n$ odd, is irreducible of highest weight $(\frac{1}{2}, \dots, \frac{1}{2})$ and accordingly, on a spin manifold, the tensor product bundle splits into two irreducible subbundles as follows:
\[\mathrm{TM}\otimes \sigm\cong\sigm\oplus \ker(c),\]
where $c:\mathrm{TM}\times \sigm\to\sigm$ denotes the \emph{Clifford multiplication} of a vector field with a spinor.
There are thus two generalized gradients: the \emph{Dirac operator} $D$, which is locally explicitly given by the formula: $D\varphi=\sum_{i=1}^{n}e_i\cdot\nabla_{e_i}\varphi, \quad \text{for all }\varphi\in\Gamma(\sigm)$ (where $\{e_i\}_{i=\overline{1,n}}$ is a local orthonormal basis and the middle dot is a simplified notation for the Clifford multiplication) and the \emph{twistor (Penrose) operator} $T$: $T_X\varphi=\nabla_X\varphi+\frac{1}{n}X\cdot D\varphi$.\\
For $n$ even, the spinor representation splits with respect to the action of the volume element into two irreducible subrepresentations, $\Sigma_n=\Sigma_n^{+}\oplus\Sigma_n^{-}$, of highest weights $(\frac{1}{2}, \dots, \frac{1}{2},\pm\frac{1}{2})$, whose elements are usually called \emph{positive}, respectively \emph{negative half-spinors}. Accordingly, there is a splitting of the spinor bundle: $\sigm=\mathrm{\Sigma^+ M} \oplus \mathrm{\Sigma^- M}$ and the decomposition of the tensor product is then given by: $\mathrm{T}^*\mathrm{M}\otimes \mathrm{\Sigma^{\pm} M}=\mathrm{\Sigma^{\mp} M}\oplus \ker(c)$.
Again the projections onto the first summand correspond to the Dirac operator and onto $\ker(c)$ to the twistor operator.
\end{Example}

\begin{Example}[Rarita-Schwinger Operator]\label{explrs}
Let $n\geq 3$ be odd and consider the so-called {\it twistor bundle}, which is the target bundle of the twistor operator acting on spinors, denoted by $\mathrm{\Sigma_{3/2} M}$. This is the vector bundle associated to the irreducible $\Spin(n)$-representation of highest weight $(\frac{3}{2},\frac{1}{2}, \cdot,\frac{1}{2})$. If $n\geq 5$, it follows from the selection rule in Lemma~\ref{selectrule} that there are four relevant weights: $0$, $-\varepsilon_1$, $+\varepsilon_1$, $+\varepsilon_2$ and the corresponding four gradient targets are: $\mathrm{\Sigma_{3/2} M}$ itself, the spinor bundle $\sigm$,  the associated vector bundles to the irreducible representations of highest weights $(\frac{5}{2}, \frac{1}{2},\dots,\frac{1}{2})$, respectively $(\frac{3}{2},\frac{3}{2},\frac{1}{2},\dots,\frac{1}{2})$. If $n=3$ the last of these targets is missing. The generalized gradient corresponding to the relevant weight $\varepsilon=0$ is denoted by $D_{3/2}:=P^{(3/2,1/2, \dots, 1/2)}_{0}$. This operator is well-known especially in the physics literature and is called the \emph{Rarita-Schwinger operator}.\\
If $n$ is even, $n=2m$, then the two bundles defined by the Cartan summand in $\mathrm{T}^*\mathrm{M}\otimes \mathrm{\Sigma^{\pm} M}$ have highest weights $(\frac{3}{2},\frac{1}{2}, \dots, \frac{1}{2}, \pm\frac{1}{2})$ and the corresponding Rarita-Schwinger operators are the generalized gradients denoted by: $D_{3/2}^{\pm}=P^{(3/2,1/2, \dots, 1/1, \pm 1/2)}_{\mp \varepsilon_m}:\Gamma(\mathrm{\Sigma_{3/2}^{\pm}M})\to\Gamma(\mathrm{\Sigma_{3/2}^{\mp}M})$.
\end{Example}

\subsection{Generalized Gradients of $G$-Structures}

Essentially the same construction as above may be used to define generalized gradients associated to a $G$-structure. On a differentiable manifold $M$ of dimension $n$, $\GL_n M$ denotes the bundle of linear frames over $M$. If $G$ is a Lie subgroup of $\GL(n, \mathbb{R})$, then a {\it $G$-structure} on $M$ is a differentiable subbundle of $\GL_n M$ with structure group $G$. The existence of a $G$-structure on a manifold is a topological condition.  As we are mainly interested in Riemannian geometry, we will consider in the sequel $G$ to be one of the subgroups of $\SO(n)$ which arise in important geometric situations and are mostly encountered in literature. These groups are exactly the ones in Berger's list of holonomy groups. Thus, in the sequel, we assume that
\begin{equation}\label{listgps}
G\in \left\{\SO(n), \mathrm{U}\left(\frac{n}{2}\right), \SU\left(\frac{n}{2}\right), \Sp\left(\frac{n}{4}\right), \Sp(1)\!\cdot\!\Sp\left(\frac{n}{4}\right), G_2, \Spin(7)\right\},
\end{equation}
where in the last two cases the dimension of the manifold is assumed to be $7$, respectively $8$. Notice that all the groups in \eqref{listgps} are compact. Moreover, their Lie algebras are simple, except for $\mathfrak{u}(\frac{n}{2})=i\mathbb{R}\oplus \mathfrak{su}(\frac{n}{2})$, which has a $1$-dimensional center and $\mathfrak{sp}(1)\oplus\mathfrak{su}(\frac{n}{4})$, \mbox{which is semisimple}.

Every finite-dimensional representation of a compact Lie group $G$ is equivalent to a unitary one, so that it can be decomposed as a direct sum of irreducible representations. Thus, without loss of generality, we consider in the sequel complex finite-dimensional irreducible representations of $G$, which are parametrized by the dominant weights. In Table~\ref{tableconfw} we wrote down the suitable positivity conditions that define the dominant weights for each group in \eqref{listgps}. The coordinates of the weights are given with respect to a chosen basis of a fixed Cartan subalgebra of the Lie algebra $\mathfrak{g}$ of $G$. 

In order to construct geometric first order differential operators of a $G$-structure, the natural starting point is, by analogy to the $\SO(n)$-case, a connection on the principal $G$-bundle. The general setting is now the following: there exists a $G$-structure on $M$, denoted by $GM$, and $\nabla$ is  a connection on $GM$; for a finite-dimensional complex irreducible $G$-representation of highest weight $\lambda$, $V_\lambda M$ denotes the associated vector bundle and $\nabla^\lambda$ the covariant derivative induced by $\nabla$ on $V_\lambda M$. We consider the restriction of the standard representation $\tau$ to the subgroup $G$ and further denote it by $\tau$. The associated vector bundle to $GM$ and $\tau$ is just the tangent bundle $\mathrm{TM}$. The real $G$-representation $\tau$ is irreducible, because the groups in \eqref{listgps} are known to act transitively on the unit sphere in $\mathbb{R}^n$. The complexification of $\tau$ remains an irreducible $G$-representation, except for the unitary and special unitary group, when it splits into two irreducible summands. 

The main ingredient needed to define the notion of $G$-generalized gradient is the following representation theoretical result (see also \cite{usgw2}):
\begin{Theorem}\label{pdecomp}
Let $G$ be one of the groups considered above, $\tau$ the restriction of the standard representation to $G$ and $\lambda$ an irreducible $G$-representation. The decomposition of the tensor product $\tau\otimes\lambda=\mathbb{R}^n\otimes_{\mathbb{R}}\lambda \cong  \mathbb{C}^n\otimes_{\mathbb{C}}\lambda$ is described as follows:
\begin{equation}\label{decmp}
\tau\otimes\lambda=\underset{\varepsilon\subset \lambda}{\oplus}(\lambda + \varepsilon),
\end{equation}
where by $\varepsilon\subset \lambda$ we denote the relevant weights of $\lambda$. Moreover, the relevant weights are described by the following selection rule: a weight $\varepsilon$ of $\tau$ is relevant for $\lambda$ if and only if $\lambda+\varepsilon$ is a dominant weight, with the exception of the weight $\varepsilon=0$, which occurs only for the groups $G_2$ and $\SO(n)$ with $n$ odd. For $\varepsilon=0$ a stronger condition must be fulfilled, namely: $\lambda-\lambda_{\tau}$, respectively $\lambda-\lambda_{\Sigma_n}$ are dominant weights, where $\lambda_\tau$ and $\lambda_{\Sigma_n}$ are the highest weight of the standard representation of $\mathfrak{g}_2$, respectively of the spinor representation of $\so(n)$. 
The decomposition \eqref{decmp} carries over to vector bundles:
\begin{equation}\label{decmpbdl}
\mathrm{T}^*\mathrm{M}\otimes_{\mathbb{R}} V_\lambda M=(\mathrm{T}^*\mathrm{M})^{\mathbb{C}}\otimes_{\mathbb{C}}V_\lambda M=\underset{\varepsilon\subset \lambda}{\oplus}V_{\lambda + \varepsilon}M.
\end{equation}
\end{Theorem}

For the proof of Theorem~\ref{pdecomp} we mention that the decomposition \eqref{decmp} is a special case of an important result in representation theory, sometimes called the general Clebsch-Gordan theorem, which provides formulas for the multiplicities of the irreducible components of the tensor product of two irreducible representations of a semisimple Lie algebra. There are different methods to compute these multiplicities, most of them being based on Weyl's character formula or, equivalently, on Kostant's formula for the multiplicity of a weight. For instance, the more particular question on which tensor products of simple Lie algebras are multiplicity-free is completely answered. The pairs of fundamental weights $\omega_1$, $\omega_2$ such that the tensor product $m_1\omega_1\otimes m_2\omega_2$ is multiplicity-free for all $m_1,m_2\geq 0$ are classified by Littelmann, \cite{ltm}, and, more generally, the multiplicity-free tensor products of simple Lie algebras have been completely classified by Stembridge, \cite{stb}, and independently, for the exceptional Lie algebras, by King and Wybourne, \cite{kw}. From their classification it follows in particular that the decomposition of $\tau\otimes\lambda$ is multiplicity-free. However, this tensor product is a very special case, because $\tau$ is the standard representation. In \cite{mpth} we give a direct proof, also based on Weyl's character formula, of the fact that the splitting \eqref{decmp} is multiplicity-free. Moreover, the argument yields the stated characterization of the relevant weights.

As the tensor product $\mathrm{T}^*\mathrm{M}\otimes V_\lambda M$ is multiplicity-free (by Theorem~\ref{pdecomp}), it follows that its decomposition is unique and the projections onto the irreducible subbundles are well-defined (otherwise one may just define the projections onto the isotypical components). This allows us to define the $G$-generalized gradients as follows: 
\begin{Definition}\label{defggg}
Let $(M,g)$ be a Riemannian manifold carrying a $G$-structure and $\lambda$ be an irreducible $G$-representation. The \emph{$G$-generalized gradient} $P^{\nabla^\lambda}_\varepsilon$, corresponding to the relevant weight $\varepsilon$ of $\lambda$, acting on sections of $V_\lambda M$, is defined by the composition:
\begin{equation}\label{eqdefggg}
\Gamma(V_\lambda M) \xrightarrow{\nabla^\lambda} \Gamma(\mathrm{TM}^*\otimes V_\lambda M) \xrightarrow{\Pi_\varepsilon} \Gamma(V_{\lambda+\varepsilon} M),
\end{equation}
where $\nabla^\lambda$ is the connection induced on $V_\lambda M$ by a $G$-connection $\nabla$ and $\Pi_\varepsilon$ is the projection onto the subbundle $V_{\lambda+\varepsilon} M$. 
\end{Definition}

The Definition~\ref{defggg} may be given more generally, for any group $G\subseteq \GL(n)$ which satisfies the technical condition \eqref{decmp}, \emph{i.e.} such that the tensor product of any $G$-irreducible representation with the restriction of the standard $\GL(n)$-representation is multiplicity-free. 

Notice that $N(\lambda):=\sharp\{\varepsilon |\, \varepsilon \text{ is a relevant weight of }  \lambda\}\leq \dim (\tau)$, so that there are at most $n$ general\-ized gradients for each dominant weight $\lambda$ and this is the generic case.

\begin{Example}[Holomorphic and Anti-Holomorphic Generalized Gradients]\label{holantihol}
If $G$ is $\mathrm{U}(\frac{n}{2})$ or $\SU(\frac{n}{2})$, then the decomposition of the complexified tangent bundle into the eigenspaces $i$ and $-i$ of the corresponding almost complex structure of the manifold: $T^{\mathbb{C}}M=T^{1,0}M\oplus T^{0,1}M$, yields a splitting of the covariant derivative: $\nabla^\lambda=(\nabla^\lambda)^{1,0}+(\nabla^\lambda)^{0,1}$. Consequently, the set of $\mathrm{U}(\frac{n}{2})$ or $\SU(\frac{n}{2})$-generalized gradients acting on sections of an irreducible vector bundle $V_\lambda M$ splits into two subsets, namely the sets of gradients factorizing over the complementary projections:
\[\xymatrix@R-25pt{
 & & & \Gamma(T^{0,1}M\otimes_\mathbb{C} V_\lambda M)\ar[r]^-{\Pi_\varepsilon} & \Gamma(V_{\lambda+\varepsilon} M) \\
\Gamma(V_\lambda M) \ar[r]^-{\nabla^\lambda} & \Gamma(T^{\mathbb{C}}M\otimes_\mathbb{C} V_\lambda M)  \ar[urr]^{\pr^{1,0}} \ar[drr]^{\pr^{0,1}} & \\
 & & & \Gamma(T^{1,0}M\otimes_\mathbb{C} V_\lambda M)\ar[r]^-{\Pi_{\varepsilon'}} & \Gamma(V_{\lambda+\varepsilon'} M),}\]
which are called {\it holomorphic}, respectively {\it anti-holomorphic} generalized gradients. This apparently skewed notation is due to the isomorphisms \mbox{$T^{1,0}\cong (T^{0,1})^*$} and $T^{0,1}\cong (T^{1,0})^*$. We notice that the weights of $T^{1,0}$ are equal to those of $T^{0,1}$ with opposite sign.\\
On a K\"ahler spin manifold of complex dimension $m$, examples of $\mathrm{U}(m)$-generalized gradients are the projections of the Dirac operator, $D^+$ and $D^-$, and the twistor operator of type $r$, $T_r$, for $r=0,\dots,m$, which act on sections of the irreducible subbundles $\s{r}$ of the spinor bundle (see \cite{mp} for details).
\end{Example}

\begin{Remark}
The exceptional case of the zero weight in the selection rule in Theorem~\ref{pdecomp} provides interesting generalized gradients, which have the same source and target bundle, so that in particular they have spectra. In this case the corresponding generalized gradient carries sections of $V_\lambda M$ to sections of a copy of $V_\lambda M$ which is a subbundle in $\mathrm{T}^*\mathrm{M}\otimes V_\lambda M$. Examples of such generalized gradients are the following: the Dirac operator, the Rarita-Schwinger operator, $*d$ acting on $\frac{n-1}{2}$-forms in odd dimension $n\geq 3$. Explicit computations of the spectra of these operators have been done on certain manifolds, for instance by Branson, \cite{br99}, on spheres. In general, there exist only estimates of the eigenvalues of these operators.
\end{Remark}

\section[Conformal Invariance]{The Conformal Invariance of $G$-Generalized Gradients}\label{confinv}

In this section we prove that the generalized gradients associated to a $G$-structure are conformally invariant. We show that the connection which ensures the conformal invariance of $G$-generalized gradients is the minimal connection of the $G$-structure. The classical conformal invariance (for the groups $\SO(n)$ and $\Spin(n)$) established by Fegan, \cite{feg}, follows then as a special case, since the minimal connection of an $\SO(n)$-structure  coincides with the Levi-Civita connection.

As mentioned above, the conformal invariance was shown by Homma, \cite{h1}, \cite{h2}, \cite{h4}, for the generalized gradients associated to $\mathrm{U}(n)$, $\Sp(n)$ and $\Sp(1)\!\cdot\!\Sp(n)$-structures. For these subgroups Homma's proof of the conformal invariance is given by explicit computations, based on the relation between the enveloping algebra of the Lie algebra of $G$ and the algebraic structure of the principal symbols of generalized gradients. Our approach uses the framework of conformal geometry and has the advantage that, on the one hand, it leads us to a more general result and, on the other hand, gives a uniform and direct proof of all known cases, avoiding the specific computations for each subgroup. 

\subsection{Main Result}
Let $(M,c)$ be an oriented $n$-dimensional manifold with a conformal structure, \emph{i.e.} an equivalence class of Riemannian metrics, where two metrics are equivalent $\bar{g}\sim g$ if there exists a function $u:M \to \mathbb{R}$ such that $\bar{g}=e^{2u}g$. In the language of $G$-structures this is equivalent to a reduction of the structure group of the tangent bundle to the conformal group $\CO(n)=\{A\in \GL(n,\mathbb{R}) | \, A^t A=a I_n, a>0 \}\cong \mathbb{R}_+^*\times \SO(n)=\{aA|\, a\in\mathbb{R}_+^*, A\in \SO(n)\}$.

Each irreducible representation of $\CO(n)$, $\tilde{\lambda}: \CO(n)\to \Aut(V)$, is identified with a couple $(\lambda,w)$, where $\lambda$ is the restriction of $\tilde{\lambda}$ to $\SO(n)$, which is still an irreducible representation, and $w$, the weight of $\tilde{\lambda}$, is determined by the restriction of $\tilde\lambda$ to $\mathbb{R}_+^*$, \mbox{which is of the form}:
\begin{equation}\label{ecw}
\tilde\lambda(a)=a^w \cdot I, \quad a\in\mathbb{R}_+^*,
\end{equation}
where $I$ is the identity of $V$ and the weight $w$ is a real or complex number, depending on whether $V$ is a real or complex representation. Let $V_{\tilde{\lambda}}M$ denote the vector bundle on $M$ associated to $\tilde{\lambda}$ and to the principal $\CO(n)$-fiber bundle of oriented $c$-orthonormal frames on $(M,c)$, denoted by $\CO_n M$, and call it of (conformal) weight $w$. Notice that any vector bundle on $M$ determined by $\GL_n M$ and a linear representation of $\GL(n,\mathbb{R})$ has a {\it natural weight}, which is the weight of the restriction of this representation to $\CO(n)$. For example, the natural conformal weight of $\mathrm{TM}$ is $1$ and of $\mathrm{T}^*\mathrm{M}$ is $-1$.

Let us first recall the definition of a Weyl structure.
\begin{Definition}
A \emph{Weyl structure} on $(M,c)$ is a linear connection $D$ on $\mathrm{TM}$, which is \emph{conformal}, \emph{i.e.} induced by a $\CO(n)$-equivariant connection on $\CO_n M$ and \emph{symmetric}, \emph{i.e.} has no torsion.
\end{Definition}

A Weyl structure $D$ is called {\it closed} if it is locally the Levi-Civita connection of a local metric in the conformal class $c$ and it is called {\it exact} if it is globally the Levi-Civita connection of a metric in $c$. The Weyl structures on the conformal manifold $(M,c)$ form an affine space modeled on the space of real $1$-forms on $M$. More precisely, two Weyl structures are related by:
\begin{equation}\label{weyldiff}
D_2=D_1+\tilde{\theta},
\end{equation}
where $\theta$ is a real $1$-form on $M$ and $\tilde{\theta}$ is the $1$-form with values in the adjoint bundle (the one associated to the adjoint representation of $\CO(n)$ on its Lie algebra $\mathfrak{co}(n)$), identified with $\theta$ by: $\tilde{\theta}(X)=\theta(X)\cdot I+\theta\wedge X$, where $\theta\wedge X$ is the skew-symmetric endomorphism defined as: $(\theta\wedge X)(Y)=\theta(Y)X-g(X,Y)\theta^{\#}$, for any metric $g$ in the conformal class and $\theta^{\#}$ the dual of $\theta$ with respect to the metric $g$.

Each conformal connection on $\mathrm{TM}$, in particular any Weyl structure $D$, induces for each linear representation $\tilde{\lambda}:\CO(n)\to \Aut(V)$ a covariant deriva\-tive $D^{\tilde{\lambda}}$ on the associated vector bundle $V_{\tilde{\lambda}} M$. If $D_1$ and $D_2$ are related by \eqref{weyldiff}, then the induced covariant derivatives $D_1^{\tilde{\lambda}}$ and $D_2^{\tilde{\lambda}}$ satisfy:
\begin{equation}\label{indweyldiff}
D_2^{\tilde{\lambda}}=D_1^{\tilde{\lambda}}+d\tilde{\lambda}(\tilde{\theta})=D_1^{\tilde{\lambda}}+\sum_{i=1}^{n} e_i^*\otimes d\lambda(\theta\wedge e_i) +w\, \theta\otimes I,
\end{equation}
where $\{e_i\}_{i=\overline{1,n}}$ is an orthonormal (conformal) frame at the point considered, $\{e_i^*\}_{i=\overline{1,n}}$ is the (algebraic) dual frame and $\tilde\lambda$ is identified as above with $(\lambda,w)$.

By $\tau$ we denote, as above, the standard representation of $\SO(n)$ on $\mathbb{R}^n$, identified with its dual $(\mathbb{R}^n)^*$, and also the representation of weight $-1$ of $\CO(n)$ on $(\mathbb{R}^n)^*$. The associated vector bundles to $\tau$ and $\CO_n M$, respectively  $\SO_g M$, are canonically identified to $\mathrm{T}^*\mathrm{M}$.

Suppose now that the Riemannian manifold $(M,g)$ admits a reduction of the orthonormal frame bundle $\SO_g M$ to a subbundle $G M$ with structure group $G$, where $G\subseteq \SO(n)$ is a closed subgroup of the special orthogonal group $\SO(n)$. First we need to establish which special connection plays for a $G$-structure the role of the Levi-Civita connection and yields the conformal invariance of generalized gradients. This is the so-called {\it minimal connection} of a $G$-structure. We now recall its definition.

Let $\mathfrak{g}\subseteq\so(n)$ be the Lie algebra of $G$ and decompose the Lie algebra $\so(n)$ of all skew-symmetric matrices as the direct sum of $\mathfrak{g}$ and its orthogonal complement: \mbox{$\so(n)=\mathfrak{g}\oplus\mathfrak{g}^\perp$}. The projections onto $\mathfrak{g}$ and $\mathfrak{g}^\perp$ are denoted by $\mathrm{pr}_{\mathfrak{g}}$, respectively $\mathrm{pr}_{\mathfrak{g}^\perp}$. The Levi-Civita connection seen as a connection form is a $1$-form $\omega^{LC}$ on $\SO_g M$ with values in the Lie algebra $\so(n)$. Restricting $\omega^{LC}$ to the subbundle $G M$ and decomposing it with respect to the above splitting, we get:
\begin{equation}\label{splitct}
\omega^{LC}=\mathrm{pr}_{\mathfrak{g}}(\omega^{LC})+\mathrm{pr}_{\mathfrak{g}^\perp}(\omega^{LC})=\omega^G+T,
\end{equation}
where $\omega^G$ is a connection form on $G M$ and $T$ is a $1$-form on $M$ with values in the associated vector bundle $G M\times_G \mathfrak{g}^\perp$. The connection corresponding to the connection $1$-form $\omega^G$ is called the {\it minimal connection} of the $G$-structure and is denoted by $\nabla^G$. $T$ is called the {\it intrinsic torsion} of the $G$-structure and is a measure for how much the $G$-structure fails to be integrable. More precisely, a $G$-structure is integrable if and only if its intrinsic torsion vanishes, which means that the Levi-Civita connection restricts to $G M$ and its holonomy group is contained in $G$. Otherwise stated, the intrinsic torsion is the obstruction for the Levi-Civita connection to be a $G$-connection. 

In the sequel we consider the generalized gradients of a $G$-structure defined by its minimal connection. They are denoted as follows:
\begin{equation}\label{defopg}
P^{G,\lambda}_{\varepsilon}=\Pi_\varepsilon\circ \nabla^{G,\lambda}: \Gamma(V_\lambda M) \to \Gamma(V_{\lambda+\varepsilon} M),
\end{equation}
where $\nabla^{G,\lambda}$ is the connection induced by the minimal connection $\nabla^G$ on $V_\lambda M$.

If $g$ and $\bar{g}$ are two conformally related metrics, there is a corresponding conformal change of the $G$-structure, denoted, with a slight abuse of notation, by $ \bar{G} M$, which is the image of $G M$ under the following principal bundle isomorphism between $\SO_g M$ and $\SO_{\bar g} M$:
\begin{equation}\label{isomgg}
\Phi^{g,\bar g}: \SO_g M \to \SO_{\bar g}M, \quad \{e_1, \dots, e_n\} \mapsto \{e^{-u}e_1, \dots, e^{-u}e_n\}.
\end{equation}
Then the following commutative diagram holds, where all the arrows are the natural inclusions:
\[\xymatrix@R-10pt{
G M\ar@{^{(}->}[r] \ar[d]& (\mathbb{R}_+^*\times G) M \ar[d]& \bar{G} M \ar@{_{(}->}[l] \ar[d]\\
\SO_g M \ar@{^{(}->}[r] & \CO_n M & \SO_{\bar{g}}M \ar@{_{(}->}[l]}\]
The right and left squares of the diagram are still commutative, when considering the minimal connections. More precisely, for instance for the left square, this means that the extension $D^G$ of the minimal connection $\nabla^G$ to $(\mathbb{R}^*_+\times G)M$ coincides with the projection onto $\mathbb{R}\oplus\mathfrak{g}\subseteq\mathfrak{co}(n)$ of the Weyl connection $D^g$ on $\CO_n M$ given by the extension of the Levi-Civita connection $\nabla^g$ to $\CO_n M$.

Any irreducible representation of $\mathbb{R}_+^*\times G$ is parametrized by $\tilde{\lambda}=(\lambda,w)$, where $\lambda$ is the highest weight of its restriction to $G$ and $w$ is its weight, determined by the restriction to $\mathbb{R}^*_+$, as in \eqref{ecw}. The tensor product $\tau\otimes\tilde{\lambda}$, where $\tau$ is the restriction of the standard representation to $\mathbb{R}_+^*\times G$, has then weight $w-1$ and it has as $\mathbb{R}_+^*\times G$ and $G$-representation the same decomposition into irreducible components of multiplicity $1$, as described by Theorem~\ref{pdecomp}.

We now consider the same construction of generalized gradients as in Definition~\ref{defggg}, starting with a Weyl structure on $\CO_n M$. Let $\lambda$ be an irreducible $G$-representation. Then, for any weight $w$, there is an irreducible representation $\tilde\lambda=(\lambda,w)$ of $\mathbb{R}_+^*\times G$ and the natural first order differential operators acting on sections of the associated bundle $V_{\tilde\lambda} M$ are defined by:
\begin{equation}\label{defopconfg}
P^{D,\tilde\lambda}_{\varepsilon}=\Pi_\varepsilon\circ D^{G,\tilde\lambda},
\end{equation}
for the connection $D^G$ given by the projection onto $(\mathbb{R}_+^*\times G) M$ of a Weyl connection $D$ on $\CO_n M$. We then have the following definition:
\begin{Definition}\label{defconfinvg}
The operator $P^{G,\lambda}_{\varepsilon}$ is called {\it conformal invariant relative to the weight $w$} if the operators defined by \eqref{defopconfg} do not depend on the Weyl structure whose projection onto the principal subbundle $(\mathbb{R}_+^*\times G) M$ defines the generalized gradients $P^{D,\tilde\lambda}_{\varepsilon}$.
\end{Definition}

It turns out that these weights, with respect to which the generalized gradients are conformally invariant, are exactly the eigenvalues of the so-called conformal weight operator. Let us first recall its definition (see \cite{feg} and \cite{usgw2}):
\begin{Definition}\label{defconfop}
The {\it conformal weight operator} of a $G$-representation $\lambda$, $\lambda:G\to \Aut(V)$, is the endomorphism defined as follows:
\begin{equation}\label{confopg}
B^{\lambda}_{\mathfrak{g}}: (\mathbb{R}^n)^*\otimes V \to (\mathbb{R}^n)^*\otimes V, \quad B^{\lambda}_{\mathfrak{g}}(\alpha\otimes v)=\sum_{i=1}^{n}e_i^*\otimes d\lambda(\mathrm{pr}_{\mathfrak{g}}(e_i\wedge\alpha))v,
\end{equation}
where $\{e_i\}_{\overline{1,n}}$ is an orthonormal basis of $\mathbb{R}^n$ and $\{e_i^*\}_{\overline{1,n}}$ the dual basis. The operator $B^{\lambda}_{\mathfrak{g}}$ is $(\tau\otimes\tilde\lambda)$-equivariant for any weight $w$, $\tilde\lambda=(\lambda,w)$. We also denote by $B^{\lambda}_{\mathfrak{g}}$ the induced endomorphism on the associated vector bundle $\mathrm{T}^*\mathrm{M}\otimes V_{\tilde\lambda} M$.
\end{Definition}

Since the algebraic endomorphism $B^{\lambda}_{\mathfrak{g}}$ is $(\mathbb{R}_+^*\times G)$-equivariant and the decomposition \eqref{decmp} is multiplicity-free, it follows from Schur's Lemma that $B^{\lambda}_{\mathfrak{g}}$ acts on each irreducible component of the decomposition \eqref{decmp} by multiplication with a scalar. 

We show in the sequel that the conformal weight operator $B^{\lambda}_{\mathfrak{g}}$ may be expressed in terms of the Casimir operators and then explicitly compute its eigenvalues. Let us recall that the {\it Casimir operator} of a $G$-representation $\lambda$ is defined by:
\begin{equation}\label{defcasop}
C^\lambda=-\sum_{\alpha}d\lambda(X_\alpha)\circ d\lambda(X_\alpha),
\end{equation}
where $\{X_\alpha\}_\alpha$ is an orthonormal basis of $\mathfrak{g}$, so that $C^\lambda$ is defined only up to a constant. In formula \eqref{bcas2} the Casimir operators are normalized with respect to the invariant scalar product induced on $\mathfrak{g}\subseteq \so(n)\cong \Lambda^2 \mathbb{R}^n$, $\langle X,Y\rangle=-\frac{1}{2}\tr(XY)$.  Usually it is convenient to compute the Casimir operators with respect to a chosen scalar product and then to renormalize them. If the representation $\lambda$ is irreducible, it follows from Schur's Lemma that $C^\lambda$ acts as a scalar $c(\lambda)$, called the {\it Casimir number}; it is a real strictly positive number, except for the trivial representation, when it is zero. The Casimir numbers may be computed by Freudenthal's formula: $c(\lambda)= \langle\lambda, \lambda+2\delta\rangle$, where $\delta$ is the Weyl vector of $\mathfrak{g}$, \emph{i.e.} is equal to half the sum of the positive roots, and are then renormalized as follows:
\[c^{\Lambda^2}(\lambda)=2\frac{\dim(\mathfrak{g})}{n}\frac{c(\lambda)}{c(\tau)}.\]

\begin{Lemma}[Fegan's Lemma, \cite{feg}] \label{feglem}
The conformal weight operator $B^{\lambda}_{\mathfrak{g}}$ is equal to
\begin{equation}\label{bcas2}
B^{\lambda}_{\mathfrak{g}}=\frac{1}{2}(C^{\tau\otimes\lambda}-I|_{\mathbb{R}^*}\otimes C^{\lambda} -C^{\tau}\otimes I|_{V}).
\end{equation}
\end{Lemma}
\noindent It is then straightforward the following:
\begin{Corollary}
The conformal weight operator $B^\lambda_{\mathfrak{g}}$ acts on each $G$-irreducible summand $\lambda+\varepsilon$ in the decomposition $\tau\otimes \lambda=\underset{\varepsilon\subset\lambda}{\oplus}(\lambda+\varepsilon)$ by multiplication with the scalar $w_\varepsilon(\lambda)$ given by:
\begin{equation}\label{cascompgen}
w_\varepsilon(\lambda)=\frac{1}{2}(c^{\Lambda^2}(\lambda+\varepsilon)-c^{\Lambda^2}(\lambda)-c^{\Lambda^2}(\tau))=\frac{\dim(\mathfrak{g})}{n}\frac{|\varepsilon|^2+2\langle\lambda+\delta,\varepsilon\rangle-\langle\tau+2\delta,\tau\rangle}{\langle\tau+2\delta,\tau\rangle}.
\end{equation}
\end{Corollary}

In particular, it follows that the eigenspaces of the conformal weight operator $B^\lambda_{\mathfrak{g}}$ are compatible with the irreducible decomposition of the tensor product $\tau\otimes\lambda$. 

We may now show the conformal invariance of $G$-generalized gradients:

\begin{Proposition}\label{confweightg}
The operator $P^{G,\lambda}_{\varepsilon}$ is conformally invariant relative to the weight $w_\varepsilon(\lambda)$ and this is the only weight with respect to which this operator is conformally invariant.
\end{Proposition}

\begin{proof}
Let $D_1$ and $D_2$ be any two Weyl connections. Then there is a real $1$-form $\theta$ on $M$ such that \eqref{weyldiff} holds. As above $D_i^G$ denotes the projection of $D_i$ onto the principal $G$-subbundle. The connections induced by $D_1^G$ and $D_2^G$ on the associated vector bundle $V_{\lambda}M$ are then related as follows:
\begin{equation}\label{indweyldiff2}
D_2^{G,\tilde{\lambda}}=D_1^{G,\tilde{\lambda}}+d\tilde{\lambda}(\mathrm{pr}_{\mathfrak{g}}(\tilde{\theta}))=D_1^{G,\tilde{\lambda}}+\sum_{i=1}^{n} e_i^*\otimes d\lambda(\mathrm{pr}_{\mathfrak{g}}(\theta\wedge e_i)) +w\, \theta\otimes I,
\end{equation}
where $\{e_i\}_{i=1,\dots,n}$ is a conformal frame and $\{e_i^*\}$ the dual frame. Thus, with respect to the conformal weight operator, we obtain:
\begin{equation}\label{weyldiffmu2}
(D_2^{G,\tilde\lambda}-D_1^{G,\tilde\lambda})(\xi)=w\, \theta\otimes\xi-B^{\lambda}_{\mathfrak{g}}(\theta\otimes\xi), \quad \text{for all }\xi\in \Gamma(V_{\tilde\lambda} M).
\end{equation}
Projecting now equation \eqref{weyldiffmu2} onto the summand $\lambda+\varepsilon$ of the decomposition of the tensor product $\tau\otimes\lambda$, we get:
\begin{equation}\label{prweyldiffl2}
(P^{D_2,\tilde\lambda}_{\varepsilon}-P^{D_1,\tilde\lambda}_{\varepsilon})(\xi)=(w-w_\varepsilon) \Pi_\varepsilon(\theta\otimes\xi), \quad \text{for all } \xi\in \Gamma(V_{\tilde\lambda} M).
\end{equation}
Hence, the generalized gradient $P^{G,\lambda}_{\varepsilon}$ is conformally invariant relative to the weight $w$ if and only if $w=w_\varepsilon(\lambda)$.\hfill $\qed$
\end{proof}

The next result expresses the conformal invariance directly in terms of the minimal connections of two conformally related $G$-structures.
\begin{Proposition}\label{equivg}
The following statements are equivalent:
\begin{enumerate}
	\item $P^{G,\lambda}_{\varepsilon}$ is conformally invariant relative to the weight $w_\varepsilon$.
  \item If $G M$ and $\bar{G} M$ are conformally related $G$-structures, for $\bar{g}=e^{2u}g$ and $G \bar{M}\hookrightarrow \SO_{\bar{g}}M$, then the corresponding generalized gradients are related by:
  \begin{equation}\label{confrelg}
  \bar{P}^{G,\lambda}_{\varepsilon}\circ \phi^{G,\bar{G}}_{w_\varepsilon}=\phi^{G,\bar{G}}_{w_\varepsilon-1}\circ P^{G,\lambda}_{\varepsilon},
  \end{equation}
where for any weight $w$, $\phi^{G,\bar{G}}_{w}$ is the isomorphism between the associated vector bundles $V^{G}_{\lambda}M:=G M\times_G V$ and $V^{\bar{G}}_{\lambda}M:=\bar{G} M\times_G V$, defined by:
  \[\phi^{G,\bar{G}}_{w}: V^{G}_{\lambda}M \to V^{\bar{G}}_{\lambda}M,\quad [(e_1, \dots, e_n), v] \mapsto [(e^{-u}e_1, \dots, e^{-u}e_n), e^{wu}v].\]
\end{enumerate}
\end{Proposition}

\begin{proof}
We consider the following diagram:
\[\xymatrix{
\Gamma(V^{G}_{\lambda}M) \ar[d]_{\nabla^{G,\lambda}} & \Gamma(V_{(\lambda,w_\varepsilon)}M)  \ar[l]_{\phi_{w_\varepsilon}^{G}}^{\sim}  \ar[d]_{D^{G,(\lambda, w_\varepsilon)}} \ar@<1ex>[d]^{D^{\bar{G},(\lambda, w_\varepsilon)}} \ar[r]^{\phi_{w_\varepsilon}^{\bar{G}}}_{\sim} & \Gamma(V^{\bar{G}}_{\lambda}M) \ar[d]^{\nabla^{\bar{G},\lambda}}\\
\Gamma(\mathrm{T}^*\mathrm{M}\otimes V^{G}_{\lambda}M) \ar[d]_{\Pi_\varepsilon}  & \Gamma(\mathrm{T}^*\mathrm{M}\otimes V_{\lambda}M)  \ar[d]^{\Pi_\varepsilon}  &  \Gamma(\mathrm{T}^*\mathrm{M}\otimes V^{\bar{G}}_{\lambda}M) \ar[d]^{\Pi_\varepsilon}\\
\Gamma(V^{G}_{\lambda+\varepsilon}M)   & \Gamma(V_{(\lambda+\varepsilon, w_\varepsilon-1)}M) \ar[l]_{\phi_{w_\varepsilon-1}^{G}}^{\sim} \ar[r]^{\phi_{w_\varepsilon-1}^{\bar{G}}}_{\sim} &  \Gamma(V^{\bar{G}}_{\lambda+\varepsilon} M)}\]
where $D^G$ is the minimal connection of the $G$-structure extended to the \mbox{$(\mathbb{R}^*_+\times G)$}-principal bundle (which may also be seen as the projection onto $\mathfrak{g}$ of the Weyl connection given by the extension of the Levi-Civita connection of the metric $g$) and $D^{G,(\lambda,w_\varepsilon)}$ is the induced connection on the vector bundle $V_{\tilde\lambda}M$ associated to the irreducible representation $\tilde\lambda=(\lambda,w_\varepsilon)$. The isomorphisms in the diagram are defined by:
\[\phi_{w}^{G}: V_{(\lambda,w)}M \to V^{G}_{\lambda}M, \quad \phi_{w}^{G}([(f_1, \dots, f_n), v])=[(e_1,\dots, e_n), a^{w}v],\]
where $f_i=a e_i$, $i=1,\dots, n$ and $\{e_i\}_{\overline{1,n}}$ is an orthonormal basis with respect to the metric $g$. With this notation we have: $\phi^{G,\bar{G}}_{w}=\phi_{w}^{\bar{G}}\circ (\phi_{w}^{G})^{-1}$. The left and right ``squares" of the diagram separately commute by the definition of the induced connection on an associated vector bundle. 

Suppose now that $(1)$ holds. Then the whole diagram commutes, since then for the weight $w_\varepsilon$ the compositions ``in the middle" give the same operator (as $P^{D,\tilde\lambda}_\varepsilon=\Pi_\varepsilon\circ D^{G,\tilde\lambda}$ does not depend on the Weyl structure $D$). The composition on its ``boundary" gives (2).

For the implication (2) $\Rightarrow$ (1) we first notice that if \eqref{confrelg} holds for a weight $w$ and for any two conformally related $G$-structures, then the above diagram is commutative. Thus, the operator $P^{D,\tilde\lambda}_{\varepsilon}$ is the same for all Weyl structures given by the minimal connections of conformally related $G$-structures. Then, \eqref{prweyldiffl2} implies that $(w-w_\varepsilon) \Pi_\varepsilon(\theta\otimes\xi)=0$, for any exact $1$-form $\theta$ on $M$ and any section $\xi\in \Gamma(V_{\tilde\lambda} M)$. At some fixed point on $M$, it follows that $(w-w_\varepsilon) \Pi_\varepsilon(\alpha\otimes v)=0$, for all $\alpha\otimes v\in(\mathrm{R}^n)^*\otimes V$, which shows that $w=w_\varepsilon$. Substituting in \eqref{prweyldiffl2}, it follows that $P^{g,\lambda}_{\varepsilon}$ is conformally invariant relative to the weight $w_\varepsilon$.\hfill $\qed$
\end{proof}

In conclusion, the explicit formula \eqref{cascompgen} allows us to compute the eigen\-values of the conformal weight operator $B^\lambda_{\mathfrak{g}}$ of a $G$-structure for any irreducible representation $\lambda$, and thus, by Propositions~\ref{confweightg} and ~\ref{equivg}, to determine the conformal weights of the $G$-generalized gradients. For completeness we give in Table~\ref{tableconfw} the explicit values of the conformal weights of $G$-generalized gradients for all subgroups $G$ in \eqref{listgps}. 

\begin{Example}[Holomorphic and Anti-Holomorphic Generalized Gradients, continued]
The eigenvalues of the conformal weight operator for $G=\mathrm{U}(\frac{n}{2})$, computed by formula \eqref{cascompgen}, are: $w_{i,-}=-\lambda_i+i-n$, $w_{i,+}=\lambda_i-i+1$, for $i=1,\dots,n$. By Proposition~\ref{confweightg}, these are the conformal weights of the corresponding $\mathrm{U}(\frac{n}{2})$-generalized gradients, the holomorphic and anti-holomorphic gradients. In \cite{h2} they are called \emph{K\"ahlerian gradients}, since if the metric is K\"ahler, the $\mathrm{U}(\frac{n}{2})$-structure is integrable and its minimal connection coincides with the Levi-Civita connection.
\end{Example}

\subsection{The Classical Conformal Invariance}
We notice that the conformal invariance in the classical case of generalized gradients defined by the Levi-Civita connection on a Riemannian (spin) manifold, which was obtained by Fegan, \cite{feg}, can be now recovered as a special case of Propositions~\ref{confweightg} and ~\ref{equivg}, considering $G=\SO(n)$ (respectively $G=\Spin(n)$).

Let $(M,g)$ be a Riemannian metric and $G=\SO(n)$. The minimal connection of the $\SO(n)$-structure is equal to the Levi-Civita connection and the corresponding generalized gradients, defined by \eqref{defgg}, are the following (here we specify the metric $g$ in the notation, because we shall compare operators associated to different metrics): 
\begin{equation}\label{defop}
P^{g,\lambda}_{\varepsilon}=\Pi_\varepsilon\circ \nabla^{g,\lambda},
\end{equation}
where $\nabla^{g,\lambda}$ is the connection induced on $V_{\lambda} M$ by the Levi-Civita connection $\nabla^g$ of $(M,g)$. We denote by $c$ the conformal class of the metric $g$ and consider the conformal manifold $(M,c)$. For any conformal weight $w$, each relevant weight $\varepsilon$ of $\lambda$ determines a family of generalized gradients, parametrized by the Weyl structures, acting on sections of the associated vector bundle $V_{\tilde\lambda} M$ (with $\tilde\lambda=(\lambda,w)$):
\begin{equation}\label{defopconf}
P^{D,\tilde\lambda}_{\varepsilon}=\Pi_{\varepsilon}\circ D^{\tilde\lambda},
\end{equation}
where $D^{\tilde\lambda}$ is the connection induced on $V_{\tilde\lambda} M$ by a Weyl connection $D$ of $(M,c)$, and $\Pi_\varepsilon$ the projection of $\mathrm{T}^*\mathrm{M}\otimes V_{(\lambda,w)} M$ onto its subbundle $V_{(\lambda+\varepsilon,w-1)} M$.

The generalized gradient $P^{g,\lambda}_{\varepsilon}$ is thus conformally invariant relative to the weight $w$ if the operators $P^{D,\tilde\lambda}_{\varepsilon}$ defined by \eqref{defopconf} do not depend on the Weyl structure $D$.

The conformal weight operator defined by \eqref{confopg} simplifies for a representation $\lambda$ of $\SO(n)$, $\lambda: \SO(n) \to \Aut(V)$, as follows:
\begin{equation}\label{confop}
B^{\lambda}: (\mathbb{R}^n)^*\otimes V \to (\mathbb{R}^n)^*\otimes V, \quad B^{\lambda}(\alpha\otimes v)=\sum_{i=1}^{n}e_i^*\otimes d\lambda(e_i\wedge\alpha)v,
\end{equation}
The Casimir operator defined by \eqref{defcasop} also simplifies for an $\SO(n)$-representation $\lambda$ as follows:
\begin{equation}\label{defcas}
C^{\lambda}=-\sum_{i<j}d\lambda(e_i\wedge e_j)\circ d\lambda(e_i\wedge e_j)
\end{equation}
and the Casimir numbers are computed by Freudenthal's formula: $c(\lambda)=\langle\lambda, \lambda+2\delta\rangle$,  where the components of $\delta$, the Weyl vector of $\so(n)$, are $\delta_i=\frac{n-2i}{2}$, for $i=1,\dots,m$. Fegan's Lemma then yields the following values of the conformal weights:
\begin{equation}\label{confw}
w_0(\lambda) =\frac{1-n}{2},\quad  w_{i,+}(\lambda):=w_{\varepsilon_i}(\lambda) =1+\lambda_i-i, \quad w_{i,-}(\lambda):=w_{-\varepsilon_i}(\lambda) =1-n-(\lambda_i-i).
\end{equation}

\begin{Remark}\label{cword}
Formulas \eqref{confw} together with the conditions on the weights to be relevant show that the conformal weights are ordered as follows:
\[w_{1,+}(\lambda)> \cdots >w_{m,+}(\lambda)>w_0(\lambda)\geq w_{m,-}(\lambda)>\cdots > w_{1,-}(\lambda),\]
if $n$ is odd, $n=2m+1$, and $w_0(\lambda)=w_{m,-}(\lambda)$ if and only if $\lambda_m=0$. However, from the selection rule given by Lemma~\ref{selectrule}, it follows that this equality case cannot occur for relevant weights, since if $\lambda_m=0$, then neither $w_0$ nor $w_{m,-}$ is a relevant weight.
For $n$ even, $n=2m$, the conformal weights are ordered as follows:
\[w_{1,+}(\lambda)>\cdots > w_{m-1,+}(\lambda)>\{w_{m,+}(\lambda), w_{m,-}(\lambda)\}> w_{m-1,-}(\lambda)>\cdots > w_{1,-}(\lambda),\]
where $w_{m,+}(\lambda)-w_{m,-}(\lambda)=2\lambda_m$, so that $w_{m,+}(\lambda)\neq w_{m,-}(\lambda)$ unless $\lambda_m=0$. Hence, the conformal weights are almost always distinct. It follows that, except for the cases of irreducible representations $\lambda$ with $\lambda_m=0$, the decomposition \eqref{decompirred} corresponds exactly to the eigenspaces of the conformal weight operator $B^{\lambda}$, which may be expanded as: \mbox{$B^\lambda=\underset{\varepsilon\subset\lambda}{\sum}w_\varepsilon(\lambda) \Pi_{\varepsilon}$}.
\end{Remark}

Propositions~\ref{confweightg} and ~\ref{equivg} yield then for $G=\SO(n)$ the following classical result (Fegan, \cite{feg}):
\begin{Proposition}\label{equiv}
Any generalized gradient $P^{g,\lambda}_{\varepsilon}$ is conformally invariant relative to the weight $w_\varepsilon$ and this is the only weight with respect to which $P^{g,\lambda}_{\varepsilon}$ is conformally invariant. With respect to two conformally related metrics $g$ and $\bar{G}$, $\bar{g}=e^{2u}g$, the conformal invariance relating the corresponding generalized gradients is expressed in the following form:
\begin{equation}\label{confrel}
P^{\bar{g},\lambda}_{\varepsilon}\circ \phi^{g,\bar{g}}_{w_\varepsilon}=\phi^{g,\bar{g}}_{w_\varepsilon-1}\circ P^{g,\lambda}_{\varepsilon},
\end{equation}
where, for any weight $w$, $\phi^{g,\bar{g}}_{w}$ is the isomorphism between the induced vector bundles defined by: $\phi^{g,\bar{g}}_{w}: V^{g}_{\lambda}M \to V^{\bar{g}}_{\lambda}M,\quad [(e_1, \dots, e_n), v] \mapsto [(e^{-u}e_1, \dots, e^{-u}e_n), e^{w u}v]$.
\end{Proposition}

The relation \eqref{confrel} expressing the conformal invariance of the generalized gradients may be rewritten in the following form, which is usually encountered in literature:
\begin{equation}\label{confrel2}
P^{\bar{g},\lambda}_{\varepsilon}=e^{(w_\varepsilon-1)u}\phi^{g,\bar{g}}\circ P^{g,\lambda}_{\varepsilon}\circ e^{-w_\varepsilon u}(\phi^{g,\bar{g}})^{-1},
\end{equation}
using the following identification of the associated vector bundles that does not use any weight, but has the advantage of being an isometry: $\phi^{g,\bar{g}}: V^{g}_{\lambda}M \to V^{\bar{g}}_{\lambda}M,\quad [s, v] \mapsto [\Phi^{g,\bar g}(s), v]$, where $\Phi^{g,\bar g}$ is the isomorphism \eqref{isomgg}.

\begin{Remark}
If $V$ is not just a representation of $\SO(n)$, but is the restriction of a representation of the general linear group $\GL(n, \mathbb{R})$, then, as noticed above, $V$ has a natural weight given by the restriction of the representation to $\CO(n)\subset \GL(n, \mathbb{R})$. In this case, one may canonically identify the associated bundles to this representation and to the principal bundles $\SO_g M$, $\SO_{\bar{g}} M$ and $\CO_n M$. When both representations of highest weight $\lambda$ and $\lambda+\varepsilon$ come from representations of $\GL(n, \mathbb{R})$, the corresponding generalized gradients associated to the conformally related metrics $g$ and $\bar g$ are related as follows (after having identified the associated vector bundles to the ones associated to $\CO_n M$):
\begin{equation}\label{confrelnat}
P^{\bar{g},\lambda}_{\varepsilon}=e^{(w_\varepsilon-\omega_{\lambda+\varepsilon}-1)u}P^{g,\lambda}_{\varepsilon}e^{-(w_\varepsilon-\omega_\lambda)u},
\end{equation}
where $w_\varepsilon$ is the conformal weight and $\omega_\lambda$, $\omega_{\lambda+\varepsilon}$ are the natural weights of $\lambda$ and $\lambda+\varepsilon$.
\end{Remark}

\begin{Example}[Generalized Gradients on Differential Forms, continued]\label{explpform}
The conformal weights of the generalized gradients corresponding to the three relevant weights of the irreducible representation $\lambda_p$
are given by \eqref{confw} as follows: $w_{p,-}(\lambda_p)=-n+p$, $w_{p+1,+}(\lambda_p)=-p$ and $w_{1,+}(\lambda_p)=1$. The relation \eqref{confrelnat} implies that the conformal invariance of these generalized gradients, acting on the vector bundles associated to $\CO_n M$, may be expressed as follows:
\[d^{\bar{g}}=d^g, \quad \delta^{\bar{g}} =e^{(-n+2p-2)u}\delta^{g}e^{(n-2p)u}.\]
The first equality just expresses the obvious fact that the exterior derivative $d$ is independent of the metric. The conformal invariance of the twistor operator is given by substituting $w_\varepsilon=w_{1,+}(\lambda_p)=1$ into \eqref{confrel2}: $T^{\bar{g}}=\phi^{g,\bar{g}}\circ T^{g}\circ e^{u}(\phi^{g,\bar{g}})^{-1}$.
\end{Example}

If $(M,g)$ is a spin manifold, the class of associated bundles is enriched by the irreducible $\Spin(n)$-representations parametrized by the dominant weights $\lambda\in(\frac{1}{2}+\mathbb{Z})^m$ and the corres\-ponding generalized gradients are also conformally invariant. The situation is very similar to the one for the group $\SO(n)$, since their Lie algebras are canonically identified. In the sequel we write down explicitly this property of conformal invariance on spin manifolds and illustrate it for the most interesting operators: the Dirac operator, the twistor (or Penrose) operator and the Rarita-Schwinger operator.

The approach is analogous to the one for the special orthogonal group, using in this case the Weyl structures on a conformal spin manifold. We recall that the \emph{conformal spin group} is the group $\mathrm{CSpin}(n)=\Spin(n)\times \mathbb{R}_+^*$, which is the universal cover of the conformal group $\CO(n)$. Its Lie algebra is $\mathfrak{cspin}_n\cong\so_n\oplus\mathbb{R}\cong\mathfrak{co}(n)$. A {\it spin structure} on a conformal manifold $(M,c)$ is given by a principal $\mathrm{CSpin}(n)$-bundle $\mathrm{CSpin}_n M$ together with a projection $\theta$, such that the following diagram commutes for every $\tilde u\in \mathrm{CSpin}_n M$:
\[\xymatrix@R-22pt{
\mathrm{CSpin}(n) \ar[r]^{a\mapsto \tilde{u}a} \ar[dd]_{\phi} & \mathrm{CSpin}_{n} M \ar[dd]^{\theta} \ar[dr]\\
& & M,\\
\CO(n) \ar[r]_{A\mapsto \theta(\tilde{u})A} & \CO_{n} M \ar[ur]}\]
where $\phi$ is the canonical projection of $\mathrm{CSpin}(n)$ onto $\CO(n)$.

If $M$ has a spin structure, then any Weyl structure $D$ on $\CO_{n} M$ induces a connection on $\mathrm{CSpin}_n M$, and therefore a covariant derivative on each associated vector bundle to a representation of $\mathrm{CSpin}(n)$. The description of the representations of $\mathrm{CSpin}(n)$ is analog to the one for $\CO(n)$: each irreducible representation $\tilde\lambda$ of $\mathrm{CSpin}(n)$ is given by a couple $(\lambda, w)$, where $\lambda$ is the restriction of $\tilde{\lambda}$ to $\Spin(n)$, which is still an irreducible representation, and $w$, the weight of $\tilde{\lambda}$, is determined by the restriction of $\tilde\lambda$ to $\mathbb{R}_+^*$: $\tilde\lambda(a)=a^w \cdot I$. Definition~\ref{defconfinvg} of conformal invariance relative to a weight may be thus carried over to generalized gradients of $\Spin(n)$. Since the formula \eqref{confop} defining the conformal weight operator $B^\lambda$ only involves the representation of the Lie algebra $\mathfrak{spin}(n)\cong \so(n)$, it is the same for the group $\Spin(n)$. Consequently, its eigenvalues, the conformal weights, are also given by \eqref{confw}.

Let $(M^n,g)$ be a Riemannian spin manifold and consider the conformal change of the metric given by $\bar{g}=e^{2u}g$. Then there is a spin structure induced on $(M^n,\bar{g})$, which is defined up to isomorphism by the following commutative diagram:
\[\xymatrix@R-10pt{
\Spin_g M \ar[r]^{\tilde{\Phi}^{g,\bar g}} \ar[d]_{\theta_g} & \Spin_{\bar{g}} M \ar[d]^{\theta_{\bar{g}}}\\
\SO_g M \ar[r]_{\Phi^{g,\bar g}} & \SO_{\bar{g}} M}.\]
Proposition~\ref{equiv} also holds for $\Spin(n)$-generalized gradients, with the only difference that for any weight $w$, the isomorphism $\phi^{g,\bar{g}}_{w}$ is replaced by the following isomorphism between vector bundles associated to $\Spin(n)$-representations:
$\tilde{\phi}^{g,\bar{g}}_{w}: V^{g}_{\lambda}M \to V^{\bar{g}}_{\lambda}M,\quad [s, v] \mapsto [\tilde{\Phi}^{g, \bar{g}}(s), e^{w u}v]$.

\begin{Example}[Dirac and Twistor Operator, continued]
The two $\Spin(n)$-generalized gradients, the Dirac operator and the twistor operator, acting on the spinor bundle $\Sigma_n M$, with $n$ odd, are conformally invariant relative to the weight $w_{0}(\rho_n)=\frac{1-n}{2}$, respectively $w_{1,+}(\rho_n)=\frac{1}{2}$, which are computed by \eqref{confw}. If $g$ and $\bar{g}$ are two conformally related metrics, $\bar{g}=e^{2u}g$, by \eqref{confrel2} we have:
\[D^{\bar{g}}\circ\tilde\phi^{g,\bar{g}}_{\frac{1-n}{2}}=\tilde\phi^{g,\bar{g}}_{-\frac{1+n}{2}}\circ D^{g},\quad
T^{\bar{g}}\circ\tilde\phi^{g,\bar{g}}_{\frac{1}{2}}=\tilde\phi^{g,\bar{g}}_{-\frac{1}{2}}\circ T^{g}.\]
In the case when $n$ is even, $n=2m$, the corresponding generalized gradients, the Dirac and the twistor operator, acting on positive and negative spinors: $D:\Gamma(\mathrm{\Sigma^{\pm} M})\to \Gamma(\mathrm{\Sigma^{\mp} M})$ and $T:\Gamma(\mathrm{\Sigma^{\pm} M})\to \Gamma(\ker(c))$ are conformally invariant relative to the conformal weights $w_{m,-}(\rho^+_n)=\!w_{m,+}(\rho^-_n)=\!\frac{1-n}{2}$, respectively  $w_{1,+}(\rho^+_n)=w_{1,+}(\rho^-_n)=\frac{1}{2}$ (computed by \eqref{confw}). 
\end{Example}

\begin{Remark}\label{rmknot}
In a simplified notation we consider:
\[\bar{\cdot}: \sigm= \Spin_g M\times_{\rho_n}\Sigma_n \to \mathrm{\Sigma \overline{M}}= \Spin_{\bar g} M\times_{\rho_n}\Sigma_n,\quad [s,\varphi]\mapsto [\tilde{\Phi}^{g,\bar{g}}(s),\varphi],\]
which is an isometry with respect to the Hermitian product on the spinor bundles. We may then rewrite the conformal invariance of $D$ and $T$ in the following more familiar expression (where $\overline D$ and $\overline T$ denote the operators associated to the metric $\bar{g}$):
\[\overline{D}(e^{-\frac{n-1}{2}u}\bar\varphi)=e^{-\frac{n+1}{2}u}\overline{D\varphi},\quad 
\overline{T}(e^{\frac{u}{2}}\bar\varphi)=e^{-\frac{u}{2}}\overline{T\varphi}.\]
We mention that the conformal invariance of $T$ is usually written in the following form: $\overline{T}_X(e^{\frac{u}{2}}\bar\varphi)=e^{\frac{u}{2}}\overline{T_X\varphi}$, where $T_X\!:\!\Gamma(\sigm)\! \to\! \Gamma(\sigm)$, so that on the right side is \mbox{the same \mbox{weight $\frac{1}{2}$}}.\\
The conformal invariance of these operators was first established by Hitchin, \cite{hit74}. The original proof is given by an explicit computation, using the following relation between the connections induced on the spinor bundles by the Levi-Civita connections of two metrics $g$ and $\bar{g}=e^{2u}g$ in the same conformal class: $\overline{\nabla}_X\overline{\varphi}=\overline{\nabla_X\varphi}-\frac{1}{2}\overline{X\cdot du\cdot\varphi}-\frac{1}{2}X(u)\overline{\varphi}$, for every $\varphi\in\Gamma(\sigm)$ and $X\in \Gamma(\mathrm{TM})$, where $\overline{X}$ is given by $\overline{X}=e^{-u}X$.
\end{Remark}

\begin{Example}[Rarita-Schwinger Operator, continued]
It follows from \eqref{confw} that the Rarita-Schwinger operator is conformally invariant relative to the same conformal weight as the Dirac operator. Namely, for $n$ odd we get: $w_0((\frac{3}{2}, \frac{1}{2}, \dots, \frac{1}{2}, \frac{1}{2}))=\frac{1-n}{2}$, and for $n$ even, $n=2m$: $w_{m,-}((\frac{3}{2}, \frac{1}{2}, \dots, \frac{1}{2}, \frac{1}{2}))=w_{m,+}((\frac{3}{2}, \frac{1}{2}, \dots, \frac{1}{2}, -\frac{1}{2}))=\frac{1-n}{2}$. Thus, the Rarita-Schwinger operator fulfills the following relation at conformal changes of the metric (with the notations in Remark~\ref{rmknot}, while here the sections $\varphi$ and $\overline\varphi$ are in the twistor bundles associated to the principal bundle $\Spin_g M$, respectively $\Spin_{\bar{g}} M$):
\begin{equation}\label{d32conf}
\overline{D_{3/2}}(e^{-\frac{n-1}{2}u}\bar\varphi)=e^{-\frac{n+1}{2}u}\overline{D_{3/2}\varphi}.
\end{equation}
\end{Example}

\subsection{Consequences}
The conformal invariance of $G$-general\-ized gradients has the following straightforward, but important consequences: 
\begin{Corollary}\label{corcinv}
Let $(M,g)$ be a Riemannian manifold admitting a $G$-structure, for a subgroup $G\subseteq \SO(n)$. Then the dimension of the kernel of any $G$-generalized gradient $P^{G,\lambda}_\varepsilon$, $\dim(\ker(P^{G,\lambda}_\varepsilon))$, is the same for all metrics conformally related to $g$.
\end{Corollary}

There are many interesting geometric objects that may be defined as sections in the kernel of $G$-generalized gradients. For instance, considering the bundle of $p$-forms and the generalized gradients in Example~\ref{pforms}, we obtain as sections in their kernel the following forms: closed forms (for $d$),  co-closed forms (for $\delta$) and conformal Killing forms (for $T$). A study of the latter, in particular examples of conformal Killing forms on nearly K\"ahler and weak $G_2$-manifolds, was given by Semmelmann, \cite{uwehabil}. If we consider the spinor bundle, as in Example~\ref{exsp}, we have as sections in the kernel of the Dirac operator the so-called \emph{harmonic spinors} and in the kernel of the twistor operator the so-called \emph{twistor spinors}. On a K\"ahler manifold such examples are provided by the so-called  \emph{K\"ahlerian twistor spinors}, which are defined as sections in the kernel of the $\mathrm{U}(\frac{n}{2})$-generalized gradient called K\"ahlerian twistor operator (see \cite{mp}). 

An important application of Corollary~\ref{corcinv} is the construction of non-trivial solutions in the kernel of $G$-generalized gradients starting from trivial ones, given for instance by parallel sections, which are well understood in terms of the holonomy representation. More precisely, any section of an associated irreducible vector bundle, which is parallel with respect to the minimal connection of the $G$-structure, induces a non-parallel section in the kernel of the $G$-generalized gradient of the conformally related $G$-structure. In the classical case, when considering the Levi-Civita connection, from parallel sections for a metric $g$ we obtain sections in the kernel of generalized gradients defined by a conformally related metric to $g$.

More generally, one may relate sections in the kernel of a generalized gradient to objects satisfying a more restrictive equation. Such an example is provided by the conformal relation between twistor and Killing spinors (see \cite{fr}), which together with the classification of manifolds admitting Killing spinors established by B\"ar, \cite{baer}, yields a description of the manifolds carrying twistor spinors. 

It turns out that the formal adjoints of generalized gradients are again generalized gradients and the following holds:
\begin{Corollary}
If $P^{G,\lambda}_\varepsilon$ is a $G$-generalized gradient with conformal weight $w_\varepsilon(\lambda)$, then its formal adjoint $P^{G,\lambda+\varepsilon}_{-\varepsilon}$ is conformally invariant with respect to the conformal weight $w_{-\varepsilon}(\lambda+\varepsilon)$.
\end{Corollary}

\begin{Remark}
We notice that in general it is not straightforward to construct higher order conformally invariant differential operators composing first order ones. The composition of two generalized gradients is not usually conformally invariant, unless the corresponding conformal weights are related by $w_{\varepsilon_2}(\lambda+\varepsilon_1)=w_{\varepsilon_1}(\lambda)-1$, in which case the composition $P^{G,\lambda+\varepsilon_1}_{\varepsilon_2}\circ P^{G,\lambda}_{\varepsilon_1}$ is a second order conformally invariant differential operator. An interesting particular case is the one when a generalized gradient is composed with its formal adjoint. For instance, the above condition is not fulfilled for the Laplace operator $\Delta$ acting on $p$-forms: it follows that $d\delta$ acting on $p$-forms is conformally invariant if and only if $p=\frac{n}{2}+1$ and similarly $\delta d$ is conformally invariant if and only if $p=\frac{n}{2}-1$, showing that $\Delta=d\delta+\delta d$ is not conformally invariant. Instead, the Laplace operator might be modified by the scalar curvature in order to make it conformally invariant. More precisely, the following formula $Y_g=4\frac{n-1}{n-2}\Delta_g+\mathrm{scal}_g$ defines the so-called \emph{conformal Laplacian} or \emph{Yamabe operator} on an $n$-dimensional Riemannian manifold $(M,g)$, for $n\geq 3$, which plays a crucial role in the solution of the Yamabe problem of finding a metric of constant scalar curvature in a given conformal class on $M$.  
\end{Remark}

We finally mention that the importance of $G$-generalized gradients also comes from the fact that they naturally give rise, by composition with their formal adjoints, to second order differential operators acting on sections of associated vector bundles. Particularly important are the extreme cases of linear combinations of such second order operators: if the linear combination provides a zero-order operator, then it is a curvature term and one obtains a so-called \wb formula (systematic approaches to the study of these formulas are provided by Homma, \cite{h3}, and by Semmelmann and Weingart, \cite{usgw2}); if the linear combination is a second order differential operator, then it is interesting to determine when it is elliptic (a complete classification of these elliptic operators for the structure groups $\SO(n)$ and $\Spin(n)$ has been given by Branson, \cite{br1}).

\newpage
\begin{landscape}
 \begin{table}[ht]
 \caption{Conformal weights of $G$-generalized gradients \vspace{0.5cm}}
 \label{tableconfw}
\scalebox{0.90}{%
 \begin{tabular}[h]{|c|c|c|c|c|}
 \hline
$\dim(M)$ & Group & Geometry & Highest weight  & Conformal Weights \\
  \hline
  \hline
 & & & & \\
\multirow{2}{*}{$2m$} & $\Spin(2m)$ & (spin) & $\lambda=(\lambda_1,\dots,\lambda_m)\in\mathbb{Z}^m \cup(\frac{1}{2}+\mathbb{Z})^m$ & $w_{i,-}=1-2m-\lambda_i+i$, $i=\overline{1,m}$ \\
 & $\SO(2m)$ &  oriented Riemannian & $\lambda_1\geq \cdots\geq \lambda_{m-1}\geq |\lambda_m|$ & $w_{i,+}=1+\lambda_i-i$, $i=\overline{1,m}$\\
   & & & & \\
  \hline
   & & & & \\
\multirow{2}{*}{$2m+1$}& $\Spin(2m+1)$ &(spin) & $\lambda=(\lambda_1, \dots, \lambda_m)\in\mathbb{Z}^m\cup(\frac{1}{2}+\mathbb{Z})^m$ & $w_0=-m$\\
 & $\SO(2m+1)$ & oriented Riemannian & $\lambda_1\geq \cdots\geq \lambda_m\geq 0$ & $w_{i,-}=-2m-\lambda_i+i$, $i=\overline{1,m}$ \\
 & & & & $w_{i,+}=1+\lambda_i-i$, $i=\overline{1,m}$ \\
   & & & & \\
  \hline
   & & & & \\
\multirow{2}{*}{$2m$} & $\SU(m)$ & (special) & $\lambda=(\lambda_1, \dots, \lambda_m)\in\mathbb{Z}^m$ & $w_{i,-}=-\lambda_i+i-m$, $i=\overline{1,m}$ \\
 & $\mathrm{U}(m)$& almost Hermitian & $\lambda_1\geq \cdots\geq \lambda_m$ & $w_{i,+}=\lambda_i-i+1$, $i=\overline{1,m}$ \\
   & & & & \\
 \hline
   & & & & \\
 \multirow{2}{*}{$4m$}&  \multirow{2}{*}{$\Sp(m)$} & almost & $\lambda=(\lambda_1, \dots, \lambda_m)\in\mathbb{Z}^m$ & $w_{i,-}=-\lambda_i+i-2m-1$, $i=\overline{1,m}$\\
& & hyper-Hermitian & $\lambda_1\geq \cdots\geq \lambda_m\geq 0$ &  $w_{i,+}=\lambda_i-i+1$, $i=\overline{1,m}$\\
   & & & & \\
 \hline
    & & & & \\
\multirow{2}{*}{$4m$} & \multirow{2}{*}{$\Sp(1)\!\cdot\!\Sp(m)$} & almost & $\beta\in\mathbb{Z}$, $\lambda=(\lambda_1, \dots, \lambda_m)\in\mathbb{Z}^m$ & $w_{\pm,i,-}=1/2(-\lambda_i\pm \beta/m +i-2m-1-1/m\pm1/m)$,\\
& & quaternion-Hermitian & $\beta\geq 0$, $\lambda_1\geq \cdots\geq \lambda_m\geq 0$  & $w_{\pm,i,+}=1/2(\lambda_i\pm \beta/m  -i+1-1/m\pm1/m)$, $i=\overline{1,m}$\\
   & & & & \\
  \hline
    & & & & \\
\multirow{3}{*}{$7$}& \multirow{3}{*}{$G_2$} & \multirow{3}{*}{$G_2$-structure} & $\lambda=(\lambda_1,\lambda_2)\in\mathbb{Z}^2$ & $w_{1,\pm}=-(5/3\mp5/3)\pm 1/3(2\lambda_1+\lambda_2)$\\
& & & $\lambda_1\geq \lambda_2 \geq 0$ & $w_{2,\pm}=-(5/3\mp 4/3)\pm 1/3(\lambda_1+2\lambda_2)$ \\
&  & & & $w_{3,\pm}=-(5/3\mp1/3)\pm 1/3(\lambda_1-\lambda_2)$\\
   & & & & \\
  \hline
     & & & & \\
\multirow{4}{*}{8} & \multirow{4}{*}{$\Spin(7)$} & \multirow{4}{*}{$\Spin(7)$-structure} &  & $w_{1,\pm}=-(9/4\mp 9/4)\pm1/2(\lambda_1+\lambda_2+\lambda_3)$\\
 & & & $\lambda=(\lambda_1,\lambda_2, \lambda_3)\in\mathbb{Z}^3$ & $w_{2,\pm}=-(9/4\mp 7/4)\pm1/2(\lambda_1+\lambda_2-\lambda_3)$\\
 & & & $\lambda_1\geq \lambda_2 \geq \lambda_3\geq 0$ & $w_{3,\pm}=-(9/4\mp3/4)\pm1/2(\lambda_1-\lambda_2+\lambda_3)$\\
 & & & & $w_{4,\pm}=-(9/4\mp1/4)\pm1/2(\lambda_1-\lambda_2-\lambda_3)$\\
 & & & & \\
\hline 
\end{tabular}}
 \end{table}
\end{landscape}

\end{document}